\documentclass[11pt]{amsart}

\usepackage{graphicx}
\usepackage{amssymb}
\usepackage{epstopdf}
\usepackage{amsmath}
\usepackage{amsthm}
\usepackage{tikz}
\usepackage{pgfplots}
\usepackage{dsfont}
\usepackage{subcaption}
\usepackage[english]{babel}
\usepackage{multirow}
\usepackage[left=2.5cm,right=2.5cm,top=2cm,bottom=2cm]{geometry}
\usepackage{hyperref}
\usepackage{wrapfig}
\usepackage{caption}
\usepackage{cleveref}
\newcommand{\R}{\mathbb{R}}
\newcommand{\e}{\varepsilon}
\newcommand{\op}{\operatorname}
\newcommand{\dx}{\op{dx}}

\newcommand{\x}{\mathrm{x}}

\crefname{equation}{Equation}{Equations}
\crefname{theorem}{Theorem}{Theorems}
\crefname{proposition}{Proposition}{Propositions}
\crefname{problem}{Problem}{Problems}
\crefname{definition}{Definition}{Definitions}
\crefname{corollary}{Corollary}{Corollaries}
\crefname{table}{Table}{Tables}
\crefname{figure}{Figure}{Figures}

\begin{document}
\title{Sensitivity analysis for dose deposition in radiotherapy via a Fokker-Planck model}
\author[R. Barnard]{Richard C. Barnard}
\address{Institute for Mathematics and Scientific Computing\\ University of Graz\\Heinrichstr 36\\8010 Graz, Austria}
\author[M. Frank]{Martin Frank}

\author[K. Krycki]{Kai Krycki}
\address{Center for Computational Engineering Sciences\\ RWTH Aachen University\\Schinkelstr. 2\\52062 Aachen, Germany}

\date{\today}
\begin{abstract}
In this paper we study the sensitivities of electron dose calculations with respect to the stopping power and the transport coefficients. We focus on the application to radiotherapy simulations. We use a Fokker-Planck approximation to the Boltzmann transport equation. Equations for the sensitivities are derived by the adjoint method. The Fokker-Planck equation and its adjoint are solved numerically in slab geometry using the spherical harmonics expansion ($P_N$) and an HLL finite volume method. Our method is verified by comparison to finite difference approximations of the sensitivities. Finally, we present numerical results of the sensitivities for the normalized average dose deposition depth with respect to the stopping power and transport coefficients, demonstrating the increasing relative sensitivities as beam energy decreases.
\end{abstract}

\maketitle

\section{Introduction}

Radiotherapy dose calculation, like many other simulations, relies on physical input data that comes from measurements or microscopic theories. In the case of dose calculations arising from electron beams, this input data includes the material properties---stopping powers, scattering cross sections, transport coefficients, etc.---of different tissues which come from various databases and have various uncertainties.  A natural question arises: what effect do these uncertainties have upon dose calculation?   To our knowledge, the influence of these uncertainties has never been studied in proper detail.  Sensitivities of dose calculation with respect to beam parameters (such as shape and energy profile) were considered in \cite{FadBalMac98, SheRog02}, however. We have reviewed some possible sources in a previous paper  \cite{OF10}. As in that paper, the cross sections we use have been extracted from ICRU libraries \cite{ICRU07}.

Our purpose here is to make a first step toward a rigorous uncertainty quantification of dose computation results under uncertain input parameters. Several methods of dose computation are available; some are Monte Carlo-based codes, such as LUKA, MCNP, EGSnrc, Penelope, and GEANT, as well as their variants XVMC, VMC++, and DPM (cf.\ \cite{SpeLew08} for a recent comparison). These have the advantage of outputting dose profiles which have close correspondence with experimental results.  However, in order to use these models for determining uncertainties, statistical methods must be used (due to the stochastic nature of the computation), potentially resulting in very large numbers of dose calculations.  The computational effort required can therefore be quite large, and possibly impractical.  Alternative methods involve using deterministic methods involving partial differential equation models \cite{Bru02,Boe98,TerKol02,HenIzaSie06}.  These models have the advantage of being computationally efficient.  We consider in this paper a somewhat simplified physical Fokker-Planck model, derived (for instance) in \cite{FHK07}, describing electrons and involving only the stopping power and the first transport coefficient.  We use this model to compute the sensitivities of the mean penetration depth (more precisely the mean penetration depth with respect to the percentage depth dose) with respect to the fully energy-dependent stopping power and transport coefficient. This is achieved via adjoint calculus, leveraging the deterministic nature of the partial differential model; sensitivities are computed via only one additional solve of the Boltzmann transport equation.  Thus, sensitivities can be obtained with significantly lower computational effort.

The rest of this paper is organized as follows.  In \cref{sec:FPeqn}, the Fokker-Planck model we use for dose calculation is described, as well as the normalized average penetration depth.  In \cref{sec:adjcalc}, the main result is presented: namely, the development of the adjoint-based calculus needed for expressing the sensitivities. This results in analytic expressions for the sensitivities of the average penetration depth with respect to the stopping power and transport coefficient.  Equipped with this, we discuss the procedure for approximating the sensitivities via the method of moments in \cref{sec:Pn}.  We then present computed sensitivities for incident beams of varying energy profiles in \cref{sec:results}. Finally, \cref{sec:discuss} consists of concluding remarks and a discussion of future questions.

\section{The Fokker-Planck Equation}\label{sec:FPeqn}
The Fokker-Planck equation treated in this paper arises as an approximate model in the dose distribution calculation in radiotherapy \cite{HJS12,OF10}. For a bounded domain $Z\subset \R^n$, it is given by
\begin{align}\label{fweq}
	-\partial_{\e}(\rho(x)S(\e)\Psi(\e, x, \Omega)) + \Omega\nabla\Psi(\e, x, \Omega) = \rho(x)T(\e)\Delta_\Omega\Psi(\e, x, \Omega),
\end{align}
where $\Psi(\e,x,\Omega)\cos(\theta)dAd\Omega d\e_e dt$ is the number of electrons moving in time $dT$ through the area $dA$ in direction $d\Omega\in S^{n-1}$ near $\Omega$  with $\theta$ the angle between $\Omega$ and the normal of $dA;$ $S(\e)$ is the stopping power; $T(\e)$ is the transport coefficient, $\rho$ is the density of the medium, and $\Delta_\Omega$ is the Laplace-Beltrami operator. 
In addition to \cref{fweq}, we impose the following conditions
\begin{align*} 
	\Psi(\infty,x,\Omega) = 0,&\qquad \text{ for } x\in Z,\Omega\in S^{n-1}\\
	\Psi(\epsilon, x,\Omega)=Q(\e,x,\Omega)&\qquad \text{ for } x\in\partial Z,~n(x)\cdot\Omega<0.
\end{align*}
 for an external source $Q(\e,x,\Omega)$ on the incoming boundary.  For numerical implementations, the ``initial condition'' that there are no particles of infinite energy is replaced with a cutoff energy.
 
   We are now interested in certain properties of the dose
\begin{align*}
	D(x) = \int_0^\infty \int_{S^{n-1}} S(\e)\Psi(\e, x, \Omega) d\Omega d\e.
\end{align*}
In particular, we want to study the dependence of the normalized average penetration depth
\begin{align*}
	\overline{\x}_\text{norm} = \frac{\int_Z x D(x)~dx}{\int_Z D(x)~dx} = \frac{\overline{x}}{\overline{D}},
\end{align*}
with respect to the stopping power and the transport coefficient. For notational simplicity, $\overline{D}$ shall denote the total dose deposition throughout the remainder of this paper and
\begin{align}\label{APD}
	\overline{\x} = \int_Z x D(x)~dx
\end{align}
the nonnormalized average penetration depth. The quantities $\overline{\x}$, $\overline{D}$, and $\overline{\x}_\text{norm}$ can be seen as functionals, mapping the energy dependent functions $S$ and $T$ to a scalar. The sensitivity of these quantities with respect to the data depends on the corresponding derivatives. 

We denote the Fr\'echet derivatives with respect to $S$ by $\overline{x}_S$, or $\partial_S \overline{x}$, and use this notation for all functionals that appear in this paper. For every fixed $S$, these objects are linear operators, acting on variations $\delta S$ in the stopping power. We denote this by 
\begin{align*}
	\overline{x}_S \big[ \delta S \big],~~\text{or } \partial_S \overline{x} \big[ \delta S \big].
\end{align*}
Finally, the chain rule for Fr\'echet derivatives yields for the quantity that we are interested in:
\begin{align*}
	\big[\overline{\x}_\text{norm}\big]_S = \frac{\overline{x}_S \overline{D} - \overline{x} \overline{D}_S}{\overline{D}^2}.
\end{align*}
We use the same notation for the derivatives with respect to $T.$  The same expressions then hold if $S$ is replaced by $T$. 

\begin{figure}
	\centering
	\begin{subfigure}{0.45\textwidth}
		\includegraphics{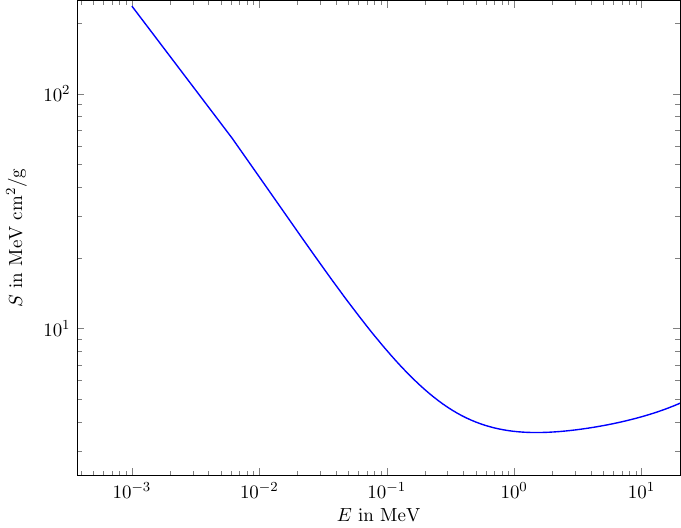}
		\caption{  Stopping Power $S$}
	\end{subfigure}\hfill
	\begin{subfigure}{0.45\textwidth}
		\includegraphics{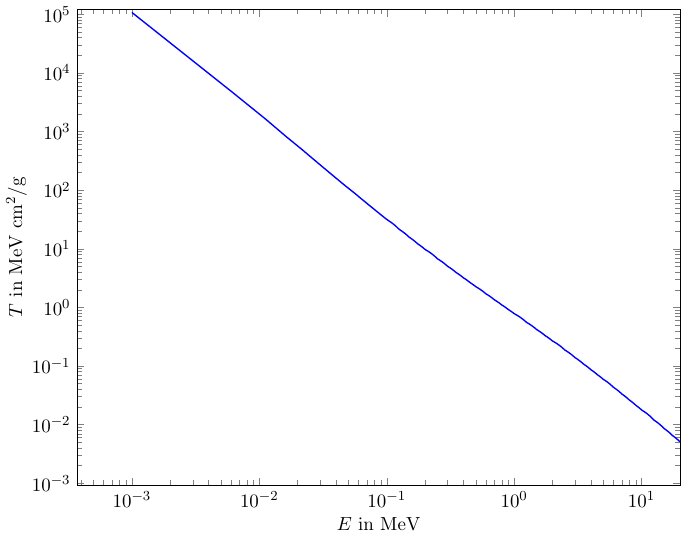}
		\caption{  Transport Coefficient $T$}
	\end{subfigure}
	\caption{Stopping Power and Transport Coefficient for Water Phantom }
	\label{fig:dosevalid}
\end{figure}

\section{Adjoint calculus and sensitivities}\label{sec:adjcalc}
In this section, we obtain analytic forms for $\big[\overline{\x}_\text{norm}\big]_T$ and $\big[\overline{\x}_\text{norm}\big]_S,$ which are, as previously mentioned, operators acting on variations in the material properties.  As such, they have representations as functions in the energy variable.  If we were instead to use a finite difference approximation, we would require $O(N_\epsilon)$ dose calculations where $N_\epsilon$ is the number of discrete energies at which we have data for $S$ and $T.$
 Additionally, inherent truncation errors in the finite difference approximation scheme tend to rise.  As we shall see in this section, adjoint methods give analytical expressions which require only a single additional solution of a Fokker-Planck equation for either the sensitivity with respect to $S$ or $T.$  

We first address the dependence of $\overline{D}$ on $S$ and $T$. Operating on Equation (\ref{fweq}) in the homogeneous case with $\rho\equiv 1$ by $\iint (\cdot) dx d\Omega$, we obtain
\begin{align*}
	\partial_\e \iint S(\e)\Psi(\e, x, \Omega)dx d\Omega = \int_{\partial Z}\int_{S^{n-1}} n\cdot\Omega \psi(\Omega,\e)d\Omega dx,
\end{align*}
as the angular integral of the Laplace-Beltrami term over the sphere vanishes.  If we assume that the incoming source is sufficiently forward peaked and the domain is large enough, we may assume that outgoing particles at the boundary are negligible.  This means we then have that 
\begin{align*}
	\partial_\e \iint S(\e)\Psi(\e, x, \Omega)dx d\Omega \approx \int_{S^{n-1}}\int_{\partial Z}\chi_{\{n\cdot\Omega<0\}} n\cdot\Omega Q(\Omega,\e)dx d\Omega .
\end{align*}
 Integrating with respect to $\e,$ we have
\begin{align*}
\overline{D}=	\iiint S(\e)\Psi(\e, x, \Omega)dx d\Omega d\e \approx \int\int_{S^{n-1}}\int_{\partial Z}\chi_{\{n\cdot\Omega<0\}} n\cdot\Omega Q(\Omega,\e)d\Omega dx d\e.
\end{align*}
for all $\e \geq 0$.  As the right hand side does not depend on $S$ or $T$, we conclude that any dependence of the total dose on $S$ or $T$ would be negligible in the case of relatively focused incoming sources. Hence, 
\begin{align*}
	\big[\overline{\x}_\text{norm}\big]_S \approx \frac{\overline{x}_S}{\overline{D}},
\end{align*}
and the main purpose of this paper is to calculate the derivatives of $\overline{x}.$

Thus, we turn here to the first order sensitivities of the average penetration depth with respect to the parameters $S$ and $T$, which are given by the derivatives of (\ref{APD}) with respect to these quantities. We define the following Lagrangian
\begin{align*}
	L:&= \overline{\x} - \langle \lambda , \Xi\rangle\\
	& = \iiint xS(\e)\Psi(\e , x,\Omega)dx d\e d\Omega - \iiint \lambda \big[ -\partial_\e (\rho S \Psi) +\Omega\nabla_x\Psi-\rho T\Delta_\Omega\Psi \big] \dx d\e d\Omega\\
	\end{align*}
where $\lambda$ is the adjoint variable which solves
\begin{align}
	-\rho S\partial_\e\lambda + \Omega\nabla_x\lambda + \rho T\Delta_\Omega\lambda + r= 0\label{eq:adjeqn}\\
	\lambda_{|\e=0} = 0,\qquad\lambda_{|n\cdot\Omega\geq 0} = 0\notag.
\end{align}
Then \begin{equation*}
L=\iiint xS\Psi -\iiint \Psi \big[ \rho S\partial_\e \lambda -\Omega\nabla_x\lambda - \rho T\Delta_\Omega\lambda\big].
\end{equation*}
Now the sensitivities of $\bar{\x}$ with respect to $S$ and $T$ can be expressed via the partial derivative of the Lagrangian with respect to these parameters. In both cases, $\lambda$ will be the solution of \cref{eq:adjeqn} with $r=xS.$  The gradient with respect to $S$, evaluated in direction $\delta S$, then reads
\begin{equation}\label{SensS}
	\partial_S L [\delta S] = \iiint [(x\Psi -\rho\Psi\partial_\e \lambda)\dx d\Omega] \delta S d\e.
\end{equation}
The gradient with respect to $T$, in direction $\delta T$  is given by
\begin{equation}\label{SensT}
	\partial_T L [\delta T] = \iiint ( \rho\Psi\Delta_\Omega\lambda \dx d\Omega)\delta T d\e.
\end{equation}

\section{Numerical method}\label{sec:Pn}
 As we are primarily concerned about the depth of dose penetration from incoming beams, we will consider \cref{fweq} and \cref{eq:adjeqn} for slab geometry. In this setting, the Laplace-Beltrami operator can be written as
\begin{align}\label{forweqSlab}
	\Delta_{\Omega} = \frac{\partial}{\partial \mu}(1-\mu^2)\frac{\partial}{\partial \mu},
\end{align}
where $\mu\in [-1,1]$ is the cosine of the polar angle from the axis of the slab geometry, and Equation (\ref{fweq}) reads
\begin{align}\label{eq:slabforw}
	-\partial_{\e}(\rho(x)S(\e)\Psi(\e, x, \mu)) + \mu\partial_x\Psi(\e, x, \mu) = \rho(x)T(\e)\frac{\partial}{\partial \mu}(1-\mu^2)\frac{\partial}{\partial \mu}\Psi(\e, x, \mu).
\end{align}
In order to reduce the dimension of the state space and develop efficient computational methods, it is common to use the method of moments in order to discretize the angular variable \cite{FraHenKla07}; for this, we multiply \cref{eq:slabforw} by the Legendre polynomials in $\mu$ and integrate with respect to $\mu.$ We define
\begin{align*}
	\Psi^{l}(\e, x):= \int_{-1}^1 P_l(\mu)\Psi(\e, x, \mu)\op{d}\mu,
\end{align*}
where $P_l$ is the $l$-th Legendre polynomial.  As the Legendre polynomials are eigenfunctions of the one-dimensional Laplace-Beltrami operator, we have after integration by parts 
\begin{align*}
	\int_{-1}^1 P_k(\mu)\frac{\partial}{\partial \mu}(1-\mu^2)\frac{\partial}{\partial \mu}\Psi(\e,x,\mu)\op{d}\mu = 
	-k(k+1)\Psi^{(k)}(\e,x)
\end{align*}
for $k = 0,\ldots,N$. 
This process results, using the well-known recursion formula for the Legendre polynomials, in a system of $N+1$ equations of the form
\begin{align}\label{PNforw}
	 -\partial_{\e}(\rho(x)S(\e)\Psi^{(k)}) + \partial_x\big(\frac{k}{2k+1}\Psi^{(k-1)} + \frac{k+1}{2k+1}\Psi^{(k+1)} \big)&= k(k+1)\rho T(\e)\Psi^{(k)},
\end{align}
where we set $\Psi^{(-1)}=0$. The system is closed by the relation $\Psi^{(N+1)} = 0$. This corresponds to a truncation of the expansion of the angular variable via a basis of Legendre polynomials. The corresponding adjoint system is
\begin{align}\label{PNadj}
	-\rho(x)S(\e)\partial_{\e}\lambda^{(0)} + \partial_x \lambda^{(1)} &= -2xS,\\
	 -(\rho(x)S(\e)\partial_{\e}\lambda^{(k)} + \partial_x\big(\frac{k}{2k+1}\lambda^{(k-1)} + \frac{k+1}{2k+1}\lambda^{(k+1)} \big)&= -k(k+1)\rho T(\e)\lambda^{(k)},
\end{align}
for $k=1,\ldots,N$, again with closure relation $\lambda^{(N+1)}=0.$  As noted in \cite{FraHerSch08}, an advantage of the $P_N$ approximation is that \cref{PNadj}, with Mark boundary conditions, is both the $P_N$ approximation of \cref{eq:adjeqn} and the adjoint arising from differentiating the semi-discretization of the Lagrangian, allowing for consistency in our discretization schemes for calculating $\partial_TL$ and $\partial_SL.$  
\subsection{Reconstructions}
Due to the non-linear nature of derivatives in (\ref{SensS}) and (\ref{SensT}), we need to approximately reconstruct the state and adjoint variables $\Psi(\e,x,\mu)$ and $\lambda(\e,x,\mu)$ in order to properly evaluate the arising integrals. Such a reconstruction can be obtained in a straightforward manner due to the linear structure of the $P_N$ expansion. As the expansion in the angular variable is assumed to be
\begin{align*}
	\tilde\Psi_{PN}(\e,x,\mu) =\sum_{l=0}^N P_l(\mu)\alpha_l (x,\e),
\end{align*}
we multiply both sides by $P_k(\mu),$ integrate, obtaining
\begin{align*}
	\alpha_l(x,\e) = \frac{1}{\langle P_l,P_l \rangle_\mu} \Psi^{(l)},
\end{align*}
with
\begin{align*}
	\langle P_l,P_l \rangle_\mu = \int_{-1}^1 P_l(\mu)P_l(\mu)\op{d}\mu = \frac{2}{2l+1}.
\end{align*}
The same reconstruction can be applied for the adjoint variable $\lambda$, too.  

\subsection{Sensitivities to Stopping Power and Transport Coefficient}
Equipped with the approximate reconstructed state and adjoints, we now turn to the approximation of \cref{SensS} and \cref{SensT}.  
We make use of the Legendre polynomials being the eigenfunctions of the Laplace-Beltrami operator. The approximate sensitivity with respect to the transport coefficient, arising from the $P_N$ approximations, simplifies \cref{SensT} to 
\begin{align*}
	\partial_T L [\delta T]_{PN} &\approx  \iiint ( \rho\tilde\Psi_{PN}\partial_\mu (1-\mu^2)\partial_\mu\tilde\lambda_{PN} \dx d\mu)[\delta T] d\e\\
	& = -\iiint (\rho \tilde\Psi_{PN}\sum_{l=0}^N\frac{2l+1}{2}l(l+1)P_l(\mu)\lambda^{(l)}d\mu \dx) [\delta T] d\e\\
	&= - \iint (\rho \sum_{l=0}^N\frac{2l+1}{2}l(l+1)\Psi^{(l)}\lambda^{(l)} \dx) [\delta T] d\e.
\end{align*}
We note that the integration with respect to $\mu$ was done analytically; the remaining integrations can be carried out numerically.

Similarly, we approximate \cref{SensS}, using the reconstructions from the $P_N$ approximation, and simplify to obtain
\begin{align*}
	\partial_S L [\delta S]_{PN} &\approx \iiint  (x\tilde\Psi_{PN} - \rho\tilde\Psi_{PN} \tilde\lambda_{PN} d\mu \dx ) [\delta S] d\e \\
	&=\iint  x\Psi^{(0)} - \rho \sum_{l=0}^N \frac{2l+1}{2} \Psi^{(l)} \partial_\e \lambda^{(l)} \dx  [\delta S] d\e.
\end{align*}
Evaluation of $\partial_{\e}\lambda^{(l)}$ and the remaining integration is done numerically.  

\section{Numerical results}\label{sec:results}
In this section we present numerical results for our method. Throughout, an explicit HLL finite volume scheme \cite{HLL} was used in order to solve the $P_{15}$ systems for both $\lambda$ and $\Psi;$ a maximum cutoff energy of  20.5 MeV was used---it is assumed no particles are present with higher energy.  A water phantom of length either 6 cm or 9 cm and density $\rho=1$ g/cm$^2$ was used.  Dose calculations for beams of form $Q_0\chi_{\mu=1}\chi_{\e=\e_c}$ for both $\e_c=5$ and $\e_c=10$ are shown, along with dose calculations using Penelope \cite{Pen09}, in \cref{fig:dosevalid}.  We see general agreement, with expected differences arising due to effects such as Bremsstrahlung not being included in the Fokker-Planck model. 

\begin{figure}
	\centering
	\begin{subfigure}{0.45\textwidth}
	\includegraphics{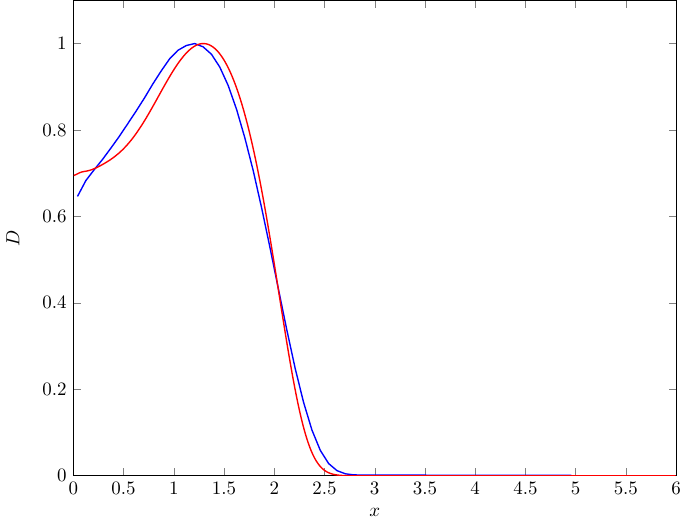}
		\caption{  Deposited dose from 5 MeV incoming beam using $P_{15}$}
	\end{subfigure}\hfill
	\begin{subfigure}{0.45\textwidth}
	\includegraphics{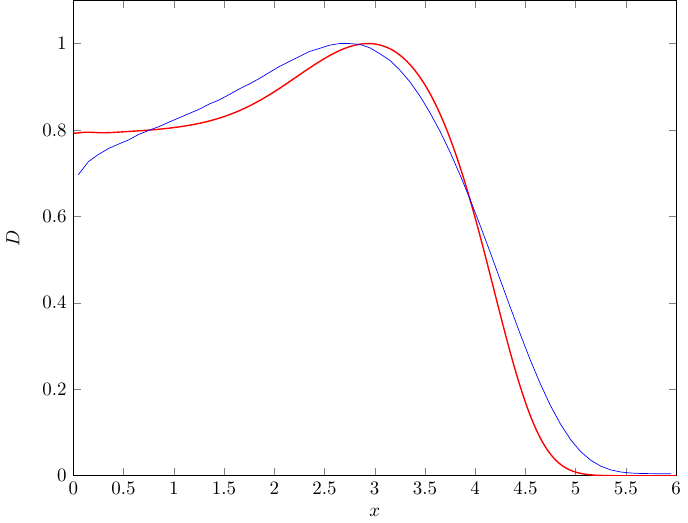}
		\caption{  Deposited dose from 10 MeV incoming beam using $P_{23}$}
	\end{subfigure}
	\caption{Percentage dose deposition (PDD) in tissue via $P_N$ approximation (red) and Penelope (blue) along axis (in cm) }
	\label{fig:dosevalid}
\end{figure}
\begin{figure}
	\centering
	\begin{subfigure}{0.45\textwidth}
\includegraphics{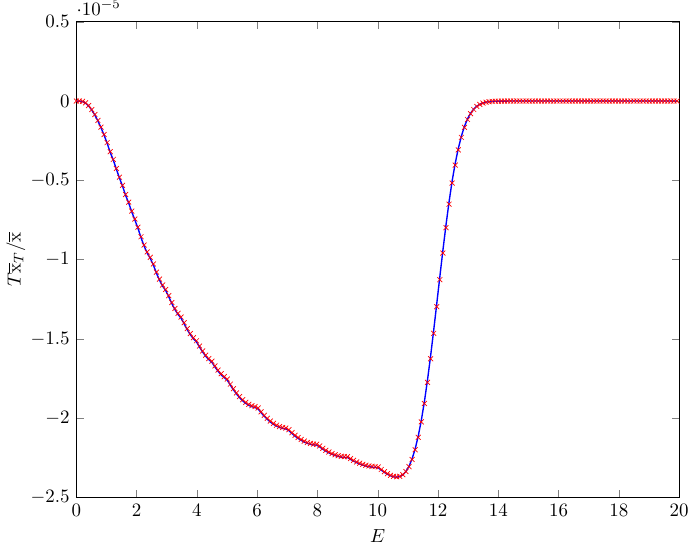}
		\caption{  $\x_T/\overline{D}$ for 10 MeV incoming beam}
	\end{subfigure}\hfill
	\begin{subfigure}{0.45\textwidth}
\includegraphics{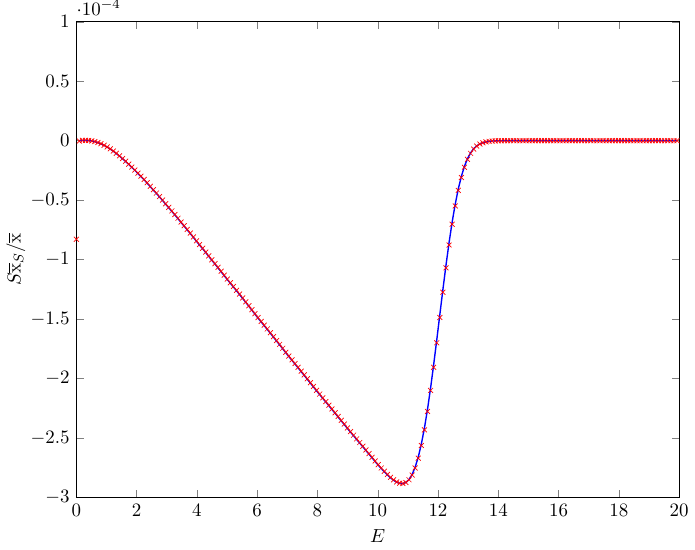}
		\caption{  $\x_S/\overline{D}$ for 10 MeV incoming beam}
	\end{subfigure}
	\caption{Derivatives of $\x$ with respect to $S$ and $T$ computed both via finite differences (in red) and via \cref{SensS} or \cref{SensT} (in blue).}
	\label{fig:10MeVSens}
\end{figure}

The relative sensitivities for the normalized penetration depth, given by
\begin{align*}
\bigg(\big[\overline{\x}_\text{norm}\big]_T\bigg)_{rel}=\frac{\overline{\x}_T}{\overline{\x}}T,\qquad \bigg(\big[\overline{\x}_\text{norm}\big]_S\bigg)_{rel}=\frac{\overline{\x}_S}{\overline{\x}}S
\end{align*} 
for a 10 MeV incoming beam are shown in \cref{fig:10MeVSens}.  For verification, a finite differences approximation was generated.  The resulting relative error in the $S$ sensitivity was $5.7697\times10^{-4}$ and in the $T$ sensitivity was $4.744\times10^{-5}$ .  
We next compute sensitivities under more realistic conditions for the incoming beam.  Beams take the form
\begin{equation*}
Q(\epsilon,\mu)= e^{\frac{-(\e-\e_c)^2}{2\sigma_{\e}^2}}e^{\frac{-(\mu-1)^2}{2\sigma_{\mu}^2}}
\end{equation*}
where $\e_c=6,9,12,16$ MeV, $\sigma_{\e}=\frac{0.1\e_c}{\sqrt{2\log2}},$ and $\sigma_{\mu}=0.1.$ The length of the slab is 6 cm with the exception of the 16 MeV beam, where 9 cm is used to account for a deeper beam penetration.  
In comparison to the sensitivities shown in \cref{fig:10MeVSens}, where there is a sharper cutoff, the effects of smearing in the energy profile of the beam are clear.  We also see a clear increase in the relative sensitivity with respect to the transport coefficient as the beam energy decreases.  The $L^\infty$ norm of the relative sensitivities for the various beams are shown in \cref{tab:varbeamdata}.


\begin{figure}
	\centering
	\begin{subfigure}{0.45\textwidth}
\includegraphics{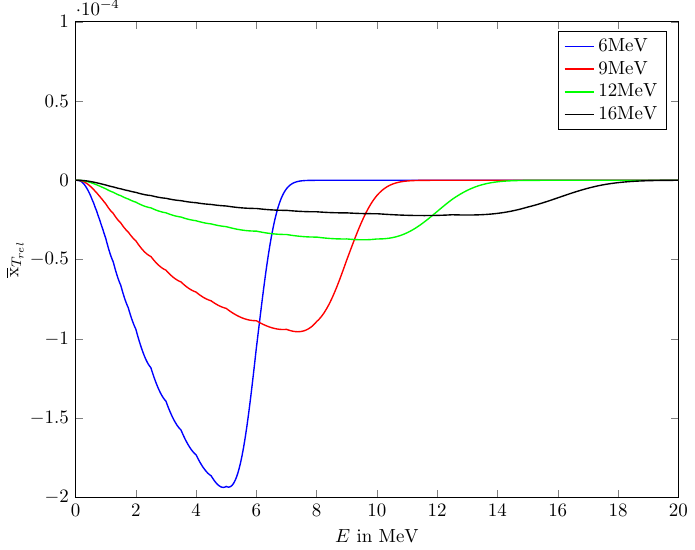}
		\caption{  Relative sensitivity of $\overline{\x}$ w.r.t. $T$ for incoming beams of varying energies}
	\end{subfigure}\hfill
	\begin{subfigure}{0.45\textwidth}
\includegraphics{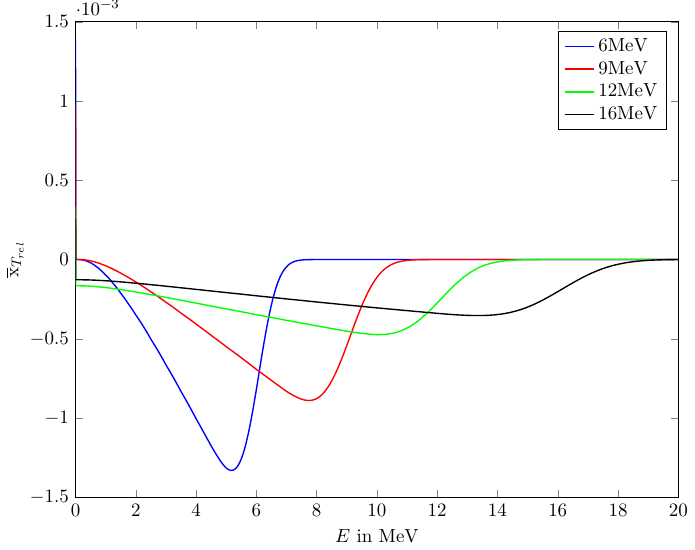}
		\caption{  Relative sensitivity of $\overline{\x}$ w.r.t. $S$ for incoming beams of varying energies}
	\end{subfigure}
	\caption{Relative sensitivities of $\overline{\x}$ with respect to $S$ and $T$ for beams centered at 6, 9, 12, and 16 MeV.}
	\label{fig:VarBeamSens}
\end{figure}

\begin{table}
\centering
 \begin{tabular}{| c || c | c | c| }
\hline
Beam energy &  $\overline{\x}/\overline{D}$ & $\|(\overline{\x}/\overline{D})_{T_{rel}}\|_\infty$ &$\|(\overline{x}/\overline{D})_{S_{rel}}\|_\infty$\\ \hline
6	&  1.2222		& $1.2964\times 10^{-4}$	&  $9.5242\times 10^{-4}$		\\ \hline
9	&  1.8316		& $6.3858\times 10^{-5}$	& $6.3463\times 10^{-4}$			\\ \hline
12 	&  2.4118		& $3.7713\times 10^{-5}$	& $4.7595\times 10^{-4}$		\\ \hline
16	&  4.7396		& $2.2259\times 10^{-5}$	& $3.5359\times 10^{-4}$		\\ \hline
\end{tabular}•
\caption{ Summary of $\overline{\x}$ and $\overline{\x}_\text{norm}$ and relative sensitivities for various beam energies.}
\label{tab:varbeamdata}
\end{table}

We see the sensitivity of the dose deposition can vary significantly in magnitude for beams of different energies.  All sensitivities reach maximum magnitude in energies slightly below the centers of the beams' energy profiles.  However, the sensitivities are also significantly stretched for the higher energy beams.  The sensitivities decrease relatively quickly once the maximum sensitivity has been reached, as fewer electrons at those energies are present in the system.  Thus, the influence of $S$ and $T$ are significantly lessened in those portions of the energy profile.  This decrease at higher energies is more rapid for beams with particles of only one energy, as seen in \cref{fig:10MeVSens}, where only trace electrons are present above 10 MeV.  In all cases, we see that the the dose penetration depth is significantly more sensitive to relative changes in the stopping power than relative changes in the transport coefficient.  This sensitivity is often of an order of magnitude higher than $\bigg(\big[\overline{\x}_\text{norm}\big]_S\bigg)_{rel}$ at its maximum.  

\section{Concluding remarks}\label{sec:discuss}
The adjoint calculus associated with the Fokker-Planck model for dose calculation led to computationally efficient methods for computing sensitivities for  the average dose penetration depth, a quantity of significant interest in radiotherapy.  This involves only two solves of a Fokker-Planck equation, which can be done efficiently via the $P_N$ method.  The sensitivities significantly vary according to both the energy profile of the incoming beam as well as its relative focus in angle.  Due to the efficiency of the computations, however, this does not pose a significant challenge, as sensitivities may be readily recomputed as new beams are investigated.  However, in all cases, the stopping power was much greater significance in computing the penetration depth.

We have so far only computed sensitivities for water phantoms in one dimension.  It would be of interest to extend this method to both problems in higher dimensions and problems involving materials with higher/lower transport coefficients such as bone and air (for computations of beams near the lung or throat).  The adjoint calculus would be unchanged; however moving to higher dimensions may require attention as the $P_N$ method may require higher numbers of moments for adequate dose computations.  Additionally, other quantities of interest (such as practical range) in the evaluation and planning of radiotherapy treatments may be studied. this would involve altering the adjoint system.  Finally, spatially inhomogeneous regions of the body would involve stopping powers and transport coefficients which are spatially dependent.  Considering this problem would lead to uncertainty quantification in the presence of uncertainties with respect to imaging and positioning of the source with respect to the patient.

\section*{Acknowledgements}
The work of Barnard  was supported in part by the Austrian Science Fund (FWF) under grant SFB F32 (SFB ``Mathematical Optimization and Applications in Biomedical Sciences'') and in part by the German Research Foundation DFG under SPP 1253 ``Optimization with partial differential equations''. 


\end{document}